\documentclass[a4paper, 11pt]{article}

\newcommand{\DRAFT}[1]{#1}
\newcommand{\FINAL}[1]{}

\usepackage[latin1]{inputenc}
\usepackage{amsmath,amssymb}
\usepackage{amsthm}
\usepackage{xspace}
\usepackage[hyphens]{url}
\usepackage{algorithmic, algorithm}
\usepackage{cite}
\usepackage[normal]{subfigure}
\usepackage{siunitx}
\usepackage{microtype}
 \usepackage[usenames,dvipsnames,svgnames]{xcolor}
\usepackage[margin=3cm]{geometry}
\usepackage{setspace}
\usepackage{empheq}

\usepackage{epstopdf}

\usepackage{multirow}
\usepackage{slashbox}
\usepackage{rotating}

\newcommand{\cut}[1]{}

\newcommand{\lcd}{}

 \newcommand{\dom}{\ensuremath{\operatorname{dom}}}
 \newcommand{\Diag}{\ensuremath{\operatorname{Diag}}}

  \newcommand{\prox}{\ensuremath{\operatorname{prox}}}
  \newcommand{\eR}{\ensuremath{\operatorname{\mathbb{R}}}}
  \newcommand{\eN}{\ensuremath{\operatorname{\mathbb{N}}}}
  \newcommand{\Cc}{\ensuremath{\operatorname{\mathcal{C}}}}
  
\newcommand{\lo}{\ensuremath{\ell_1}\xspace} 
\newcommand{\lt}{\ensuremath{\ell_2}\xspace} 
\newcommand{\loot}{\ensuremath{\ell_1/\ell_2}\xspace} 
\newcommand{\sot}{\ensuremath{\text{SOOT}}\xspace} 

\newtheorem{proposition}{Proposition}

\newtheorem{definition}{Definition}

\title{Euclid in a Taxicab:  Sparse Blind Deconvolution with  Smoothed \loot Regularization} 

\author{Audrey~Repetti\footnote{Universit\'e Paris-Est, LIGM, CNRS-UMR 8049.}, Mai~Quyen~Pham$^{*,\dagger}$, Laurent~Duval\footnote{IFP Energies nouvelles.}, \\ \'Emilie~Chouzenoux$^*$, and Jean-Christophe~Pesquet$^*$
}

\date{}

\begin{document}
\maketitle
\begin{abstract}
The \loot ratio regularization function has shown good performance for retrieving sparse signals in a number of recent works, in the context of blind deconvolution. Indeed, it benefits from a scale invariance property much desirable in the blind context. 
However, the \loot function raises some difficulties when solving the nonconvex and nonsmooth minimization
problems resulting from the use of such a penalty term in current restoration methods.
In this paper, we propose a new penalty based on a smooth approximation to the \loot function. In addition, we develop a proximal-based algorithm to solve variational problems involving this function and \lcd{we derive theoretical convergence results}. We demonstrate the effectiveness of our method through a comparison with a recent alternating optimization strategy dealing with the exact \loot term, on an application to seismic data blind deconvolution.   
\end{abstract}
%


\pagestyle{myheadings}
\thispagestyle{plain}
\markboth{ }{ }

\section{Introduction}
\label{sec:intro}

\begin{figure}[tb]
\centering
\includegraphics[width=8.6cm]{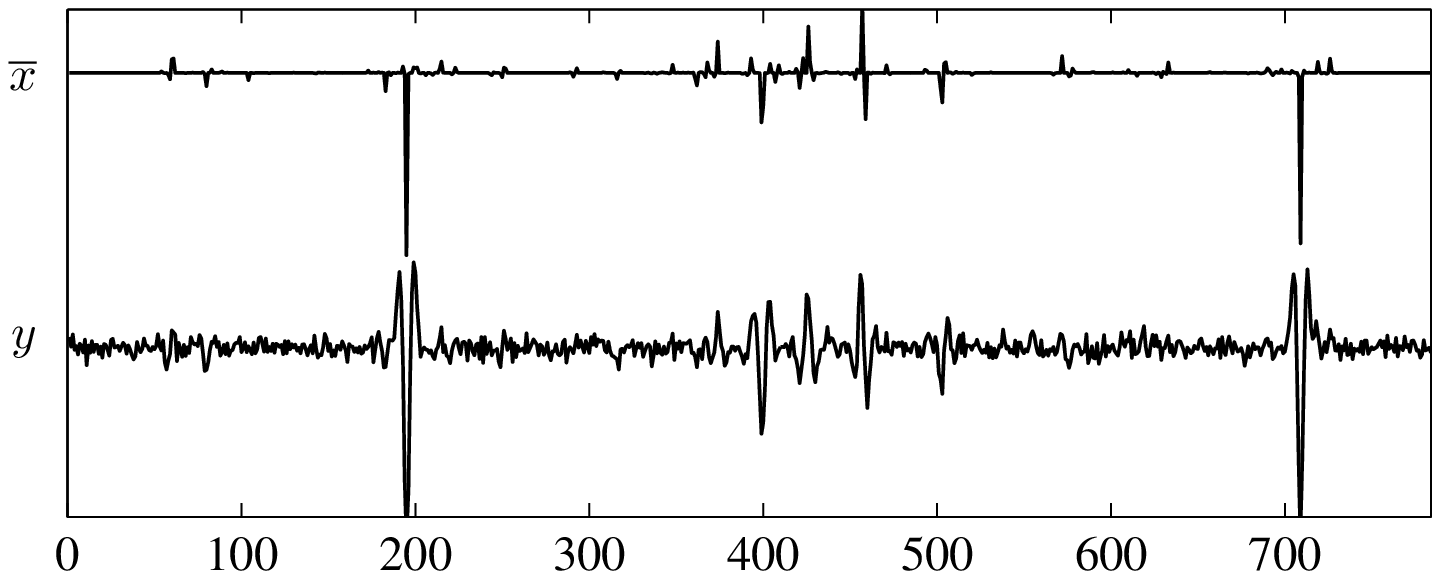} 
\caption{Unknown  seismic signal  $\bar{x}$ (top), blurred/noisy observation $y$ (bottom). \label{fig:observation}}
\end{figure}

Many experimental settings are modeled as inverse problems. They resort  to estimating an unknown signal $\overline{x} \in \eR^N$ from observations $y \in \eR^N$,  
through the measurement process:
\begin{equation}	\label{eq:conv}
y = \overline{h} \ast \overline{x} + w\,,
\end{equation}
an illustration of which is provided in Fig.~\ref{fig:observation}.
Here, $\overline{h} \in \eR^S$ represents an impulse response (e.g. a linear sensor response  or a  ``blur'' convolutive point spread function), $\ast$ denotes a discrete-time convolution operator (with appropriate boundary processing), and $w \in \eR^N$ is a realization of a random variable modeling an additive noise. Standard approaches, such as Wiener filtering and its statistical extensions \cite{Pesquet_J_2009_j-ieee-tsp_sure_adsidp},  
aim at minimizing criteria based on the squared Euclidean norm ($\ell_2^2$). 
However, the use of the sole least squares data fidelity term is prone to noise sensitivity and 
the addition of an $\ell_2^2$ regularization often leads to over-smoothed estimates. 
The deconvolution problem becomes blind, even more ill-posed, when the blur kernel $\overline{h}$  is  unknown, and needs to be estimated as well as the target signal. Applications include communications (equalization or channel estimation) \cite{Haykin_S_1994_book_blind_d}, 
nondestructive testing \cite{Nandi_A_1997_j-ieee-tsp_blind_dusnta}, 
geophysics 
\cite{Kaaresen_K_1998_j-geophysics_multichannel_bdss,Takahata_A_2012_j-ieee-spm_unsupervised_pgsrskabdbss,Pham_M_2014_j-ieee-tsp_primal-dual_proximal_astbafasmr}, image processing \cite{Kundur_D_1996_j-ieee-spm_blind_id,Kundur_D_1996_j-ieee-spm_blind_idr,Kato_M_1999_j-ieice-tfeccs_set-theoretic_bidbhsdm,Ahmed_A_2014_j-ieee-tit_blind_ducp}, medical imaging and remote sensing \cite{Campisi_P_2007_book_blind_idta}. 
Blind deconvolution, being an underdetermined problem, often requires additional hypotheses. A usual approach  seeks  estimates $(\widehat{x} , \widehat{h}) \in \eR^N \times \eR^S$ of $(\overline{x},\overline{h} )$ as minimizers of 
the sum of a data fidelity term and additional regularization terms on the signal and on the blur kernel. Such regularization functions  account for  a priori assumptions one  imposes on original sought objects, like sparsity,  and ensure the stability of the solution.
Blind deconvolution is subject to scaling ambiguity, and suggests scale-invariant contrast functions \cite{Comon_1996_j-ieee-spl_contrasts_mbd,Moreau_E_1997_j-ieee-spl_generlaized_cmbdls}. 

A decade ago, a Taxicab-Euclidean norm ratio (\loot) arose as a sparseness measure \cite{Zibulevsky_2001_j-neural-comput_blind_sssdsd,Hoyer_P_2004_j-mach-learn-res_non-negative_mfsc,Hurley_N_2009_j-ieee-tit_comparing_ms,Barak_B_2014_p-stoc_rounding_sofsr}, used in NMF (Non-negative Matrix Factorization) \cite{Morup_M_2008_p-iscas_approximate_l0cnnmtf}.  
Earlier mentions of a one-norm/two-norm ratio   deconvolution appeared in geophysics \cite{Gray_W_1978_tr_variable_nd}. It has since  been used to constrain sharp images through wavelet frame coefficients \cite{Ji_H_2012_j-acha_image_ducsiwd}, or for sparse recovery \cite{Demanet_L_2014_j-information-inference_scaling_lrses}. 
Such a regularization term  is moreover suggested in \cite{Benichoux_A_2013_p-icassp_fundamental_pbdssip} to avoid common pitfalls in blind
sparse deconvolution.

Recently, \cite{Krishnan_D_2011_p-cvpr_blind_dnsm} proposed an alternating minimization algorithm to deal with the \loot regularization function. Its originality consists of transforming the \loot nonconvex regularization term   into a convex \lo regularization function. This is done in a reweighted fashion, by fixing the denominator \lt from the previous iterate. An iterative shrinkage-thresholding algorithm finally solves the remaining \lo regularized problem. 
Although the convergence of this approach has not been deeply investigated, it appears to be quite efficient in practice. More recently, \cite{Esser_E_2013_j-siam-j-imaging-sci_method_fsssnnlspa} proposed a scaled gradient projection algorithm for minimizing a smooth approximation of the \loot function, however limited to the case when the sparse signal to retrieve takes nonnegative values. 
We generalize this idea to a parametrized Smoothed One-Over-Two (\sot) penalty for signed, real data. 
We present 
a novel efficient method based on recent results in nonconvex optimization combining an alternating minimization strategy with a forward-backward iteration \cite{Bolte_J_2013_j-math-prog-ser-a_proximal_almnnp,Chouzenoux_E_2013_tr_block_cvmfba}. 
Moreover, we accelerate the convergence of our algorithm by using a Majorize-Minimize (MM) approach  \cite{Sotthivirat_S_2002_j-ieee-tip_image_rupspscaa,Chouzenoux_E_2013_tr_block_cvmfba,Chouzenoux_E_2014_j-optim-theory-appl_variable_mfbamsdfcf}. 
Section~\ref{sec:Problem} introduces the minimization problem. Section~\ref{sec:Algo} describes the proposed method and provides 
 convergence results.
The algorithm performance, compared  with \cite{Krishnan_D_2011_p-cvpr_blind_dnsm}, is discussed  in Section~\ref{sec:Appli} for seismic data blind deconvolution. 
Some conclusions are  drawn in Section~\ref{sec:Conclusion}.

\section{Optimization model}
\label{sec:Problem}

\subsection{Optimization tools}

Our minimization strategy relies on two optimization principles. 
Let $U\in\eR^{M \times M}$ be a symmetric positive definite (SPD) matrix. Firstly, we define the $U$-weighted proximity operator \cite[Sec. XV.4]{HiriartUrruty_J_1993_book_convex_ama}, \cite{Combettes_P_2013_j-non-linear-anal_variable_mqfm} of a proper, lower semicontinuous, convex function $\psi \colon \eR^M \to ]-\infty,+\infty]$ at $z\in \eR^M$, relative to the metric induced by $U$, and denoted by $\prox_{U,\psi}(z)$, as the unique minimizer of $ \psi+ \frac{1}{2} \| \cdot -z \|^2_U  $, where  $\| . \|_U$ denotes the weighted Euclidean norm, i.e., 
$
(\forall z \in \eR^M)$ $\| z \|_U =  \left(z^\top U z \right)^{1/2}$. 
When $U$ is equal to $\operatorname{I}_M$, the identity matrix of $\eR^{M \times M}$, then 
$\prox_{\operatorname{I}_M, \psi}$
reduces to the original definition of the proximity operator in \cite{Moreau_J_1965_j-bull-soc-math-fr_proximite_deh}. We refer to \cite{Chaux_C_2007_j-inv-prob_variational_ffbip,Pustelnik_N_2011_j-ieee-tip_parallel_pairuhr,Combettes_P_2011_incoll_proximal_smsp} for  additional details on proximity operators.
Secondly, we introduce the Majoration-Minimization (MM) principle:
\begin{definition}	\label{def:maj}
Let $\zeta \colon \eR^M \to \eR$ be a differentiable function. Let $z \in \eR^M$. Let us define, for every $z' \in \eR^M$,
\begin{equation*}
q(z',z) = \zeta(z) + (z'-z)^\top \nabla \zeta (z) + \dfrac{1}{2} \| z'-z\|_{U(z)}^2,
\end{equation*}
where $U(z)\in \eR^{M \times M}$ is a Semidefinite Positive (SDP) matrix. 
Then, $U(z)$ satisfies the majoration condition for $\zeta$ at $z$ if $q(\cdot, z)$ is a quadratic majorant of the function $\zeta$ at $z$, i.e., for every $z' \in \eR^M$, $\zeta(z') \le q(z',z)$.
\end{definition}
If function $\zeta$ has an $L$-Lipschitzian gradient on a convex subset $C \subset \eR^M$, with $L>0$, i.e., for every $(z,z') \in C^2 $, 
$\| \nabla \zeta(z) - \nabla \zeta(z') \| \le L \| z-z' \|$, 
then, for every $z \in C$, a quadratic majorant of $\zeta$ at $z$ is trivially obtained by taking $U(z) = L \operatorname{I}_M$.

\subsection{Proposed criterion}
\label{ssec:MinimizationPb}

From now on, definitions and properties apply for every $x = (x_n)_{1\le n \le N}\in \eR^N $ and  $h \in  \eR^S$, unless otherwise stated. We propose to define an estimate $(\widehat{x} , \widehat{h})$ of $(\overline{x},\overline{h})$ as 
a minimizer of the following penalized criterion: 
\begin{equation}
F(x,h) = \rho(x,h)  +  g(x,h) + \varphi(x),\label{Pb:min}
\end{equation}
where $\rho(x,h) = \frac{1}{2} \| h \ast x -  y \|^2$ is the least-squares objective function,
$g$ introduces additional a priori information on the sought objects, and $\varphi$ models the One-Over-Two norm ratio non-convex penalty function \cite{Slavakis_K_2013_j-ieee-tsp_generalized_tosalus}, 
defined as the quotient of 
$\ell_1(x) = \sum_{n=1}^N |x_n|$ and $\ell_2(x) = \left(\sum_{n=1}^N x_n^2\right)^{1/2}$.
The resulting regularization term is both nonconvex and nonsmooth, so that finding a minimizer of $F$ is a challenging task.

The smooth approximations of \lo and \lt,  $\ell_{1,\alpha}$ (sometimes called hybrid \lo-\lt or hyperbolic penalty) and $\ell_{2,\eta}$,  are defined as follows with parametric  constants  $(\alpha, \eta)$:
\[
\ell_{1,\alpha}(x) = \sum_{n=1}^N \left( \sqrt{x_n^2 + \alpha^2} - \alpha \right),	\,
\ell_{2,\eta}(x) = \sqrt{\sum_{n=1}^N x_n^2 + \eta^2}.
\]
Note that \lo and \lt are recovered 
for $\alpha = \eta = 0$. 
We thus propose to replace the nonsmooth function \loot by 
a manageable smooth approximation.
More precisely, we employ the following surrogate function:
\begin{equation}
\varphi(x) = \lambda \log \left( \dfrac{\ell_{1,\alpha}(x) + \beta}{\ell_{2,\eta}(x)} \right),
\label{eq:reg_sparse} 
\end{equation}
with $ (\lambda,\beta,\alpha,\eta) \in ]0,+\infty[^4$. 

The $\log$ function both makes the penalty easier to handle 
and, through its concavity, tends to strengthen the sparsity promoting effect of the \loot function.
$F$ corresponds to the Lagrangian function associated with the minimization of 
$\rho  + g$ 
under the constraint
\begin{equation}	\label{eq:ConstraintS'}
 \log \left(\frac{\ell_{1,\alpha}(x)+\beta}{\ell_{2,\eta}(x)} \right) \le \log(\vartheta),
\end{equation}
for some positive constant $\vartheta $. 
Owing 
to the monotonicity of the $\log$ function, \eqref{eq:ConstraintS'} is equivalent to $(\ell_{1,\alpha}(x)+\beta) / \ell_{2,\eta}(x) \le \vartheta$, 
which, according to \eqref{eq:reg_sparse}, can be interpreted as a smooth approximation of an \loot upper bound constraint, for $\beta$ small enough. Finally, remark that lengthy but straightforward calculations allowed us to prove that $\varphi$ has a Lipschitzian gradient on any bounded convex subset of $\eR^N$, which is a desirable property for deriving an efficient algorithm to minimize~\eqref{Pb:min}. 
In the following, 
we  assume that $g$ can be split as  
\begin{equation}	\label{eq:g}
g(x,h) = g_1(x) + g_2(h),
\end{equation}
where $g_1$ and $g_2$ are (non-necessarily smooth) proper, lower semicontinuous, convex functions, continuous on their domain.
Moreover, we  denote by 
\begin{equation}	\label{eq:Lform}
f(x,h) = \rho(x,h) + \varphi(x), 
\end{equation}
the smooth part of the criterion, and 
$\nabla_1 f(x,h) \in \eR^N$ (resp. $\nabla_2 f(x,h) \in \eR^S$) the partial gradient of $f$ with respect to the variable $x$ (resp. $ h$) computed at $ (x,h) $.

\section{Proposed alternating optimization method}
\label{sec:Algo}

\subsection{Proposed \sot algorithm}

To minimize \eqref{Pb:min}, one can exploit the block-variable structure of 
$F$ by using an alternating forward-backward algorithm \cite{Luo_Z_1992_j-optim-theory-appl_convergence_cdmcdm,Bolte_J_2010_p-icip_alternating_pabir,Xu_Y_2013_j-siam-j-imaging-sci_block_cdmrmoantfc,Bolte_J_2013_j-math-prog-ser-a_proximal_almnnp,Chouzenoux_E_2013_tr_block_cvmfba}. 
At each iteration $k \in \eN$, this algorithm updates $x^k$ (resp. $h^k$) with a gradient step on $f(\cdot,h^k)$ (resp. $f(x^k,\cdot)$) followed by a proximity step on $g_1$ (resp. $g_2$). 

We use this alternating minimization method combined with an MM strategy, as described in \cite{Chouzenoux_E_2013_tr_block_cvmfba}. 
For every $(x,h)\in \eR^N \times \eR^S$, let us assume the existence of SPD matrices $A_1(x,h) \in \eR^{N \times N}$ and  $A_2(x,h) \in \eR^{S \times S}$ such that $A_1(x,h)$ (resp. $A_2(x,h)$) satisfies the majoration condition for $f(\cdot,h)$ at $x$ (resp. $f(x,\cdot)$ at $h$). Then, the \sot algorithm for the minimization of \eqref{Pb:min} is described  in Algorithm~\ref{algo:BCVMFB}.
\begin{algorithm}
\caption{\sot algorithm. \label{algo:BCVMFB}}
\begin{algorithmic}
\STATE  For every $k \in \eN$, let $J_k \in \eN^*$, $I_k \in \eN^*$ and let $(\gamma_x^{k,j})_{0 \le j \le J_k-1}$ and $ (\gamma_h^{k,i})_{0 \le i \le I_k-1}$ be positive sequences. Initialize with  $x^0 \in \dom g_1$ and $h^0 \in \dom g_2$.
\STATE \textbf{Iterations:} 
\STATE
$	
\begin{array}{l}
\text{For } k = 0,1,\ldots	\\
\left\lfloor               
\begin{array}{l} 
	x^{k,0} = x^k , \, \, h^{k,0} = h^k , \\
	\text{For } j = 0, \ldots, J_k-1\\
	\left\lfloor
	\begin{array}{l} 
		\widetilde{x}^{k,j} =  x^{k,j} - \gamma_x^{k,j} A_1(x^{k,j},h^k)^{-1} \nabla_1 f ( x^{k,j}, h^k ) , \\
		x^{k,j+1}
			= \prox_{(\gamma_x^{k,j})^{-1} A_1(x^{k,j},h^k), g_1} \left( \widetilde{x}^{k,j} \right),
	\end{array}
	\right.   \\
	x^{k+1} = x^{k,J_k}.	\\
	\text{For } i = 0, \ldots, I_k-1\\
	\left\lfloor
	\begin{array}{l} 
		\widetilde{h}^{k,i} = h^{k,i} - \gamma_h^{k,i} A_2(x^{k+1},h^{k,i})^{-1} \nabla_2 f ( x^{k+1}, h^{k,i} ),\\
		h^{k,i+1}
			= \prox_{(\gamma_h^{k,i})^{-1} A_2(x^{k+1},h^{k,i}), g_2} \left( \widetilde{h}^{k,i} \right), 
	\end{array} 
	\right.   \\
	h^{k+1} = h^{k,I_k}.	
\end{array} 
\right.             
\end{array}
$
\end{algorithmic}
\end{algorithm}
%
Note that PALM algorithm \cite{Bolte_J_2013_j-math-prog-ser-a_proximal_almnnp} is recovered as a special case if $J_k \equiv I_k \equiv 1$ and, at each iteration,  
the Lipschitz constant of $\nabla_1 f(\cdot, h^k)$ (resp. $\nabla_2 f( x^{k+1}, \cdot)$) is substituted for $A_1(x^{k,0},h^k)$ (resp. $A_2(x^{k+1},h^{k,0})$). 
However, recent works on variable metric strategies \cite{Chouzenoux_E_2013_tr_block_cvmfba,Chouzenoux_E_2014_j-optim-theory-appl_variable_mfbamsdfcf} show that the use of more judicious preconditioning matrices can significantly accelerate  the convergence of the algorithm. 
An example of such matrices is proposed in Section~\ref{ssec:MMstrat}.
Moreover, we  show in our experimental part the practical interest in terms of convergence speed of taking the number of inner loops $(I_k)_{k \in \eN}$ or $(J_k)_{k \in \eN}$ greater than one.

The convergence of Algorithm~\ref{algo:BCVMFB} can be derived from the general results established in~\cite{Chouzenoux_E_2013_tr_block_cvmfba}:

\begin{proposition}	\label{th:convergence}
Let $(x^k)_{k \in \eN}$ and $(h^k)_{k \in \eN}$ be sequences generated by Algorithm~\ref{algo:BCVMFB}. Assume that:
\begin{enumerate}

\item	\label{ass:i3}
There exists $(\underline{\nu}, \overline{\nu}) \in ]0, +\infty[^2$ such that, for all $k \in \eN$,
\begin{align*}
(\forall j \in \{0, \ldots, J_k-1\})\quad
&	\underline{\nu} \operatorname{I}_N \preceq A_1(x^{k,j},h^k) \preceq \overline{\nu} \operatorname{I}_N ,	\\
(\forall i \in \{0, \ldots, I_k-1\})\quad
&	\underline{\nu} \operatorname{I}_S \preceq A_2(x^{k+1},h^{k,i}) \preceq \overline{\nu} \operatorname{I}_S.
\end{align*}

\item	\label{ass:ii}
Step-sizes $(\gamma_x^{k,j})_{k\in\eN,0\le j \le J_k-1}$ and $(\gamma_h^{k,i})_{k\in \eN,0\le i \le I_k-1}$ are chosen in the interval $[\underline{\gamma},2-\overline{\gamma}]$
where $\underline{\gamma}$ and $\overline{\gamma}$ are some given positive real constants.

\item	\label{ass:iii}
$g$ is a semi-algebraic function.\footnote{Semi-algebraicity is a property satisfied by a wide class of functions,
which means that their graph is a finite union of sets defined by a finite number of polynomial inequalities. In particular, it is satisfied for the \sot penalty, for standard numerical implementations of the $\log$ function.}
\end{enumerate}
Then, the sequence $(x^k, h^k)_{k \in \eN}$ converges to a critical point $(\widehat{x}, \widehat{h})$ of \eqref{Pb:min}. 
Moreover, $\big(F(x^k,h^k)\big)_{k \in \eN}$ is a nonincreasing sequence converging to $F(\widehat{x}, \widehat{h})$.
\end{proposition}

\subsection{Construction of the quadratic majorants}
\label{ssec:MMstrat}

The numerical efficiency of the \sot algorithm relies on the use of quadratic majorants providing tight approximations to the criterion and whose curvature matrices are simple to compute. 
The following proposition allows us to propose SDP
matrices $A_1$ and $A_2$ for building majorizing approximations of $f$ with respect to $x$ and $h$. 
\begin{proposition}
For every $(x,h) \in \eR^N \times \eR^S$, let
\begin{align*}
& A_1(x,h) = \left( L_1(h) + \frac{9\lambda}{8\eta^2} \right) \operatorname{I}_N + \frac{\lambda}{\ell_{1,\alpha}(x)+\beta} A_{\ell_{1,\alpha}}(x),	\\
& A_2(x,h) = L_2(x)  \operatorname{I}_S ,
\end{align*}
where 
\begin{equation}	\label{eq:maj_l1s}
	A_{\ell_{1,\alpha}}(x) = \Diag \left( \left( (x_n^2 + \alpha^2)^{-1/2} \right)_{1 \le n \le N} \right) ,
\end{equation} 
and $L_1(h)$ (resp. $L_2(x)$) is a Lipschitz constant for
 $\nabla_1 \rho( \cdot,h)$ (resp. $\nabla_2 \rho(x, \cdot)$).\footnote{Such Lipschitz constants are straightforward to derive since
$\rho$ is a quadratic cost.} 
Then,  $A_1(x,h)$ (resp. $A_2(x,h)$) satisfies the majoration condition for $f(\cdot,h)$ at $x$ (resp. $f(x,\cdot)$ at $h$).
\end{proposition}
\proof{
Let us decompose $\varphi= \varphi_1 + \varphi_2$ with
$\varphi_1(x) = \lambda \log \left( \ell_{1,\alpha}(x) + \beta \right)$ and $\varphi_2(x) = - \lambda \log \left( \ell_{2,\eta}(x) \right)$.
It then suffices to prove that,
for every $x \in \eR^N$,
\begin{enumerate}
\item[(i)]	
$ A_{\varphi_1}(x) = \frac{\lambda}{\ell_{1,\alpha}(x)+\beta} A_{\ell_{1,\alpha}}(x) $ 
satisfies the majoration condition for $\varphi_1$ at $x$,

\item[(ii)]	
$\varphi_2$ has a $\mu$-Lipschitzian gradient, with $ \mu = \frac{9\lambda}{8 \eta^2}$.
\end{enumerate}

On the one hand, setting
$\tau(x) = \ell_{1,\alpha}(x) + \beta$, we have \cite{Allain_M_2006_j-ieee-tip_global_lchqa}
\begin{equation}	\label{pr:maj_ell}
\tau(x') 
\le \tau(x) + (x'-x)^\top \nabla \tau(x) + \frac{1}{2} \| x' - x \|_{A_{\ell_{1,\alpha}}(x)}^2,
\end{equation}
for every $x' \in \eR^N$, where $A_{\ell_{1,\alpha}}(x)$ is given by \eqref{eq:maj_l1s}. 

On the other hand, for every $(u,v) \in ]0, +\infty[^2$,
\begin{equation}
\log v \le \log u + \dfrac{v}{u}-1 = \log u + \dfrac{v-u}{u}. \label{eq:ineqlog}
\end{equation}
By taking $v = \tau(x') > 0$ and $u = \tau(x) > 0$, 
and by combining \eqref{pr:maj_ell} and \eqref{eq:ineqlog}, we obtain 
\DRAFT{
\begin{equation*}
\varphi_1(x') \le \varphi_1(x) + \frac{\lambda}{\tau(x)} (x'-x)^\top \nabla \tau(x) 
+  \frac{1}{2} (x'-x)^\top \frac{\lambda}{\tau(x)} A_{\ell_{1,\alpha}}(x) (x'-x).
\end{equation*}
}
\FINAL{
\begin{multline*}
\varphi_1(x') \le \varphi_1(x) + \frac{\lambda}{\tau(x)} (x'-x)^\top \nabla \tau(x) \\
+  \frac{1}{2} (x'-x)^\top \frac{\lambda}{\tau(x)} A_{\ell_{1,\alpha}}(x) (x'-x).
\end{multline*}
}
Thus, Statement (i) is proved by remarking that $\nabla \varphi_1(x) = \frac{\lambda}{\tau(x)} \nabla \tau(x)$ and $A_{\varphi_1}(x) = \frac{\lambda}{\tau(x)} A_{\ell_{1,\alpha}}(x)$.
On the other hand, the Hessian of $\varphi_2$ is given by
\begin{equation*}
\nabla^2 \varphi_2(x) = \frac{2 \lambda}{\ell^4_{2,\eta}(x)}  xx^\top - \frac{\lambda}{\ell^2_{2,\eta}(x)} \operatorname{I}_N.
\end{equation*}
Noting that 
$\ell^2_{2,\eta}(x) = \|x\|^2 +\eta^2$, 
and applying the triangular inequality yield
\begin{equation*}
\| \nabla^2 \varphi_2(x) \|
\le 	\frac{2 \lambda \|x\|^2}{( \| x \|^2 + \eta^2 )^2} + \frac{\lambda}{\| x \|^2 + \eta^2} = \chi(\|x\|),
\end{equation*}
where $\chi \colon u \in [0,+\infty[ \mapsto  \lambda \frac{3u^2+\eta^2}{(u^2 + \eta^2)^2}$. 
The derivative of $\chi$ is given, for every $u \in [0,+\infty[$, by
\begin{equation*}
\dot{\chi}(u) = \lambda \frac{2u}{(u^2+\eta^2)^3} ( \eta^2 - 3u^2 ),
\end{equation*}
thus $\chi$ is an increasing function on $[0, \eta/ \sqrt{3}]$ and a decreasing function on 
$]\eta/\sqrt{3}, + \infty[$, 
and $\sup_{u \in [0,+\infty[} \chi(u) = \chi\left( \eta / \sqrt{3} \right) =  \frac{9\lambda}{8\eta^2}$.
Hence, the proof of Statement (ii).
~$\blacksquare$
}

\section{Application to seismic data deconvolution}
\label{sec:Appli}

\subsection{Problem statement}
As some of the earliest  mentions of   \loot   deconvolution appeared in geophysics \cite{Gray_W_1978_tr_variable_nd}, blind seismic  deconvolution (or inversion \cite{Osman_O_1996_book_seismic_ssem,Ulrych_T_2005_book_information-based_ipa}) is  a natural application.
The sparse seismic signal $\overline{x} $, of length $N=784$, on the top of Fig.~\ref{fig:observation} is composed of a sequence of spikes termed primary reflection coefficients \cite{Walden_A_1986_j-geophys-prospect_nature_ngprcsd}. 
This reflectivity series indicates, in reflection seismology at normal incidence, the  travel time of seismic waves between two seismic reflectors, and the amplitude of the seismic events reflected back to the sensor. The observed seismic trace $y$ displayed in Fig.~\ref{fig:observation}-bottom follows Model \eqref{eq:conv}. 
In this context, the blur $\overline{h} $ is 
related to the generated  seismic source. We use here  a band-pass ``Ricker'' seismic wavelet (or Mexican hat \cite{Ricker_N_1940_j-geophysics_form_nswss}) of size $S = 41$ (Fig.~\ref{fig:estimation}-bottom)  with a frequency spectrum concentrated between \num{10} and \SI{40}{\hertz}.  The additive noise $w$ is a realization of a zero-mean white Gaussian noise with variance $\sigma^2$. 
Since the reflectivity series  is sparse, but limited in amplitude, we choose $g_1$ as the indicator function of the convex hypercube $[x_{\min}, x_{\max}]^N$. 
Similarly, as the seismic wavelet possesses finite energy, $g_2$ is equal to the indicator function of the set $\Cc = \{ h \in [ h_{\min}, h_{\max} ]^S \, | \, \| h \| \le \delta \}$, where $\delta>0$, and $h_{\min}$ (resp. $h_{\max}$) is the minimum (resp. maximum) value of $\overline{h}$.

\subsection{Numerical results}

\begin{figure}[bht]
\begin{center}
\includegraphics[width=8.6cm]{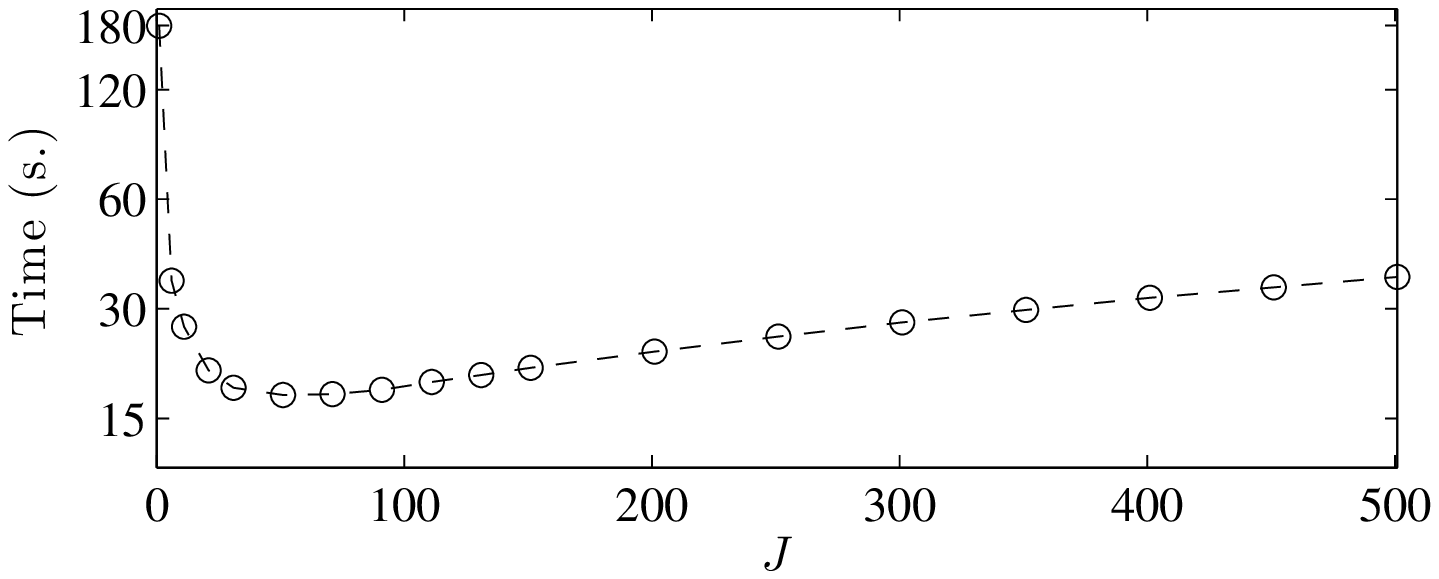} 
\caption{Reconstruction time for different numbers of inner-loops $J_k \equiv J$ (average over thirty noise realizations). \label{fig:block}}
\end{center}
\end{figure}

Fig.~\ref{fig:block} presents the variations of the reconstruction time, in seconds, with respect to the number of inner-loops $J_k \equiv J$, with $I_k \equiv 1 $ and noise level $\sigma = 0.03$. 
The reconstruction time corresponds to the stopping criterion $ \|x^k - x^{k-1}\| \le \sqrt{N}\times\num{e-6}$.
One can observe that the best compromise in terms of convergence speed is obtained for an intermediate number of inner-loops, namely $J = 71$. 
Note that the quality of the reconstruction is stable for each choice of $J$.

We gather comparisons of  the \sot algorithm with \cite{Krishnan_D_2011_p-cvpr_blind_dnsm} in Table~\ref{tab:snrxh}, where the same initialization strategy has been used for both algorithms: $x^0$ is a constant-valued signal such that $\|x^0\| \le \max \{ |x_{\min}|, |x_{\max}| \}$,
and $h^0$ is a centered Gaussian filter, such that $h^0 \in \mathcal{C}$. 
Results presented in this table, for each noise level $\sigma$, are averaged over two hundred noise realizations. 
The regularization parameters of \cite{Krishnan_D_2011_p-cvpr_blind_dnsm} and $(\lambda,\alpha,\beta,\eta)\in ]0,+\infty[^4$ of \eqref{eq:reg_sparse} are adjusted so as to minimize 
the \lo norm between the original  and the reconstructed signals. We also set, for every $k\in\eN$, $J_k=71$ and $I_k=1$. 
If both methods yield tremendous improvements in \lt and \lo norms, the \sot algorithm exhibits better results, for all noise levels, for both $\overline{x}$ and $\overline{h}$ estimates, especially in terms of $\ell_1$ norm. 
 Interestingly, the \sot algorithm is also significantly faster in this application. 

\begin{figure}[h!]
\centering
\begin{tabular}{@{}c@{}}
\includegraphics[width=8.6cm]{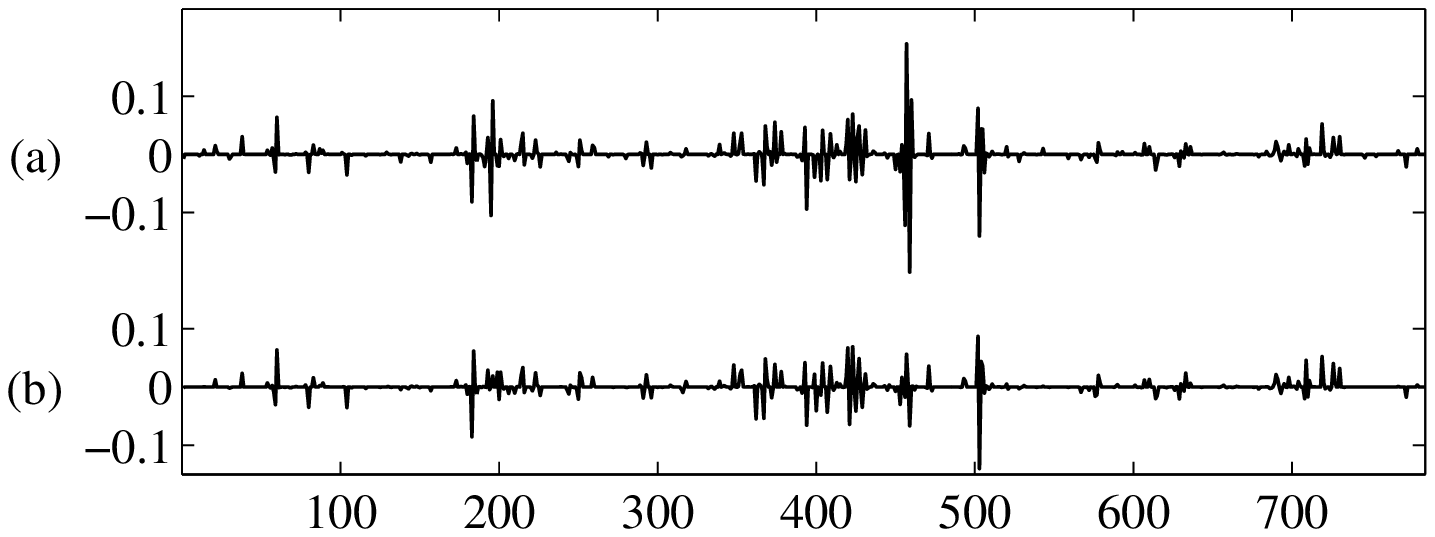} 
\\
\includegraphics[width=8.6cm]{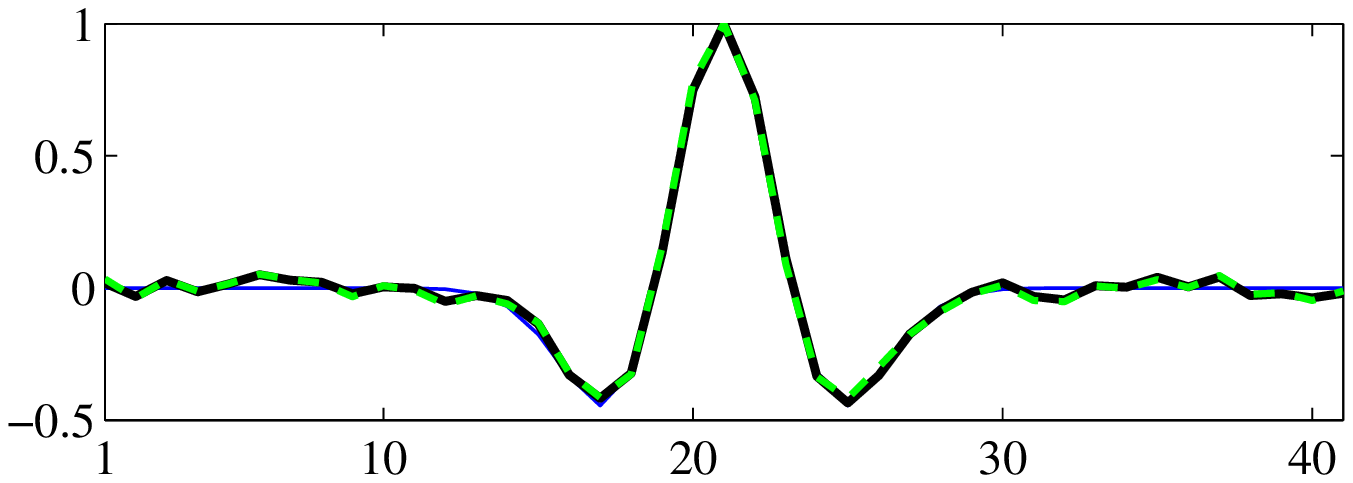} 
\end{tabular}
\caption{Top:  signal estimation error $\overline{x}- \widehat{x} $ with estimates  $\widehat{x}$  given by  \cite{Krishnan_D_2011_p-cvpr_blind_dnsm} (a) and \sot (b). Bottom: Original blur $\overline{h}$ (continuous thin blue), estimated $\widehat{h}$ with \sot (continuous thick black) and \cite{Krishnan_D_2011_p-cvpr_blind_dnsm} (dashed thick green).}
\label{fig:estimation}
\end{figure}

The performance is further assessed by subjective results for  $\sigma = 0.03$. Fig.~\ref{fig:estimation}-top shows the residual error of the sparse signal estimation $\overline{x} - \widehat{x}$, for a given noise realization, where $\widehat{x}$ is estimated with \cite{Krishnan_D_2011_p-cvpr_blind_dnsm} in (a), and with \sot in (b). 
It appears, in this example, that the error is smaller using \sot algorithm. 
The estimated blur kernels
look similar for both methods, 
as displayed in Fig.~\ref{fig:estimation}-bottom.

\begin{table}[htb]
\centering
\small{
\begin{tabular}{|c|c|c|c|c|c|}
\hline
\multicolumn{3}{|c|}{Noise level ($\sigma$)} & 0.01 & 0.02 & 0.03 \\
\hline
\hline
 \multicolumn{2}{|c|}{\multirow{2}{*}{Observation error}} 
 & \lt ($\times \num{e-2}$) & 7.14  &  7.35 & 7.68 \\
\cline{3-6}
 \multicolumn{2}{|c|}{}
 & \lo ($\times \num{e-2}$) & 2.85  &  3.44 & 4.09 \\
\hline
\hline
 \multirow{4}{*}{Signal error}
 &\multirow{2}{*}{\cite{Krishnan_D_2011_p-cvpr_blind_dnsm}}
 & \lt ($\times \num{e-2}$) & 1.23 & 1.66 & 1.84 \\
\cline{3-6}
 && \lo ($\times \num{e-3}$) & 3.79  & 4.69 & 5.30 \\
\cline{2-6}
 &\multirow{2}{*}{\sot} 
 & \lt ($\times \num{e-2}$) & 1.09  & 1.63 & 1.83\\
\cline{3-6}
 && \lo ($\times \num{e-3}$) & 3.42  & 4.30 & 4.85\\
\hline
\hline
 \multirow{4}{*}{Kernel error}
 &\multirow{2}{*}{\cite{Krishnan_D_2011_p-cvpr_blind_dnsm}}
 & \lt ($\times \num{e-2}$) & 1.88  & 2.51 & 3.21 \\
\cline{3-6}
 && \lo ($\times \num{e-2}$) & 1.44  & 1.96 & 2.53 \\
\cline{2-6}
 &\multirow{2}{*}{\sot} 
 & \lt ($\times \num{e-2}$) & 1.62  & 2.26 & 2.93 \\
\cline{3-6}
 && \lo ($\times \num{e-2}$) & 1.22  & 1.77 & 2.31 \\
\hline
\hline
 \multirow{2}{*}{ Time (s.)} 
 &\multicolumn{2}{|c|}{\cite{Krishnan_D_2011_p-cvpr_blind_dnsm}} 
 & 106 & 61 & 56 \\
 \cline{2-6}
 &\multicolumn{2}{|c|}{\sot}
 &  56  & 22 &  18 \\
\hline
\end{tabular}
}
\caption{ Comparison between   \cite{Krishnan_D_2011_p-cvpr_blind_dnsm} and \sot for  $\overline{x}$ and $\overline{h}$ estimates 
(Intel(R) Xeon(R) CPU E5-2609 v2@2.5GHz 
using Matlab~8). 
 \label{tab:snrxh}}
\vspace{-0.8cm}
\end{table}

\section{Conclusion}
\label{sec:Conclusion}

The proposed \sot for minimizing an \loot penalized criterion has been demonstrated to be quite effective in  a blind   deconvolution application on  seismic reflectivity data. 
In addition, one of its advantages is that it offers theoretically guaranteed convergence.
In future works, its use should be investigated for a broader class of application areas, where norm ratios are beneficial: adaptive filtering \cite{Loganathan_P_2009_j-ieee-taslp_class_scaec},  compression \cite{Drugman_T_2014_j-ieee-spl_maximum_pmslps}, sparse system identification \cite{Yukawa_M_2013_p-iswcs_sparsity-based_drpapfbsa}, sparse  recovery \cite{Demanet_L_2014_j-information-inference_scaling_lrses},  or cardinality-penalized clustering \cite{Chang_X_2014_PREPRINT_sparse_kmlilophddc}. 
The application of the method using a nonquadratic data fidelity term, in association with more sophisticated preconditioning matrices, is also of main interest.

\section*{Acknowledgment}

L. Duval would like to thank Igor Carron (Nuit Blanche) for useful discussions and references.


\begin{thebibliography}{10}

\bibitem{Pesquet_J_2009_j-ieee-tsp_sure_adsidp}
J.-C. Pesquet, A.~Benazza-Benyahia, and C.~Chaux,
\newblock ``A {SURE} approach for digital signal/image deconvolution
  problems,''
\newblock {\em IEEE Trans. Signal Process.}, vol. 57, no. 12, pp. 4616--4632,
  Dec. 2009.

\bibitem{Haykin_S_1994_book_blind_d}
S.~Haykin, Ed.,
\newblock {\em Blind Deconvolution},
\newblock Prentice Hall, 1994.

\bibitem{Nandi_A_1997_j-ieee-tsp_blind_dusnta}
A.~K. Nandi, D.~Mampel, and B.~Roscher,
\newblock ``Blind deconvolution of ultrasonic signals in nondestructive testing
  applications,''
\newblock {\em IEEE Trans. Signal Process.}, vol. 45, no. 5, pp. 1382--1390,
  1997.

\bibitem{Kaaresen_K_1998_j-geophysics_multichannel_bdss}
K.~F. Kaaresen and T.~Taxt,
\newblock ``Multichannel blind deconvolution of seismic signals,''
\newblock {\em Geophysics}, vol. 63, no. 6, pp. 2093--2107, Nov. 1998.

\bibitem{Takahata_A_2012_j-ieee-spm_unsupervised_pgsrskabdbss}
A.~K. Takahata, E.~Z. Nadalin, R.~Ferrari, L.~T. Duarte, R.~Suyama, R.~R.
  Lopes, J.~M.~T. Romano, and M.~Tygel,
\newblock ``Unsupervised processing of geophysical signals: A review of some
  key aspects of blind deconvolution and blind source separation,''
\newblock {\em IEEE Signal Process. Mag.}, vol. 29, no. 4, pp. 27--35, Jul.
  2012.

\bibitem{Pham_M_2014_j-ieee-tsp_primal-dual_proximal_astbafasmr}
M.~Q. Pham, L.~Duval, C.~Chaux, and J.-C. Pesquet,
\newblock ``A primal-dual proximal algorithm for sparse template-based adaptive
  filtering: Application to seismic multiple removal,''
\newblock {\em IEEE Trans. Signal Process.}, vol. 62, no. 16, pp. 4256--4269,
  Aug. 2014.

\bibitem{Kundur_D_1996_j-ieee-spm_blind_id}
D.~Kundur and D.~Hatzinakos,
\newblock ``Blind image deconvolution,''
\newblock {\em IEEE Signal Process. Mag.}, vol. 13, no. 3, pp. 43--64, May
  1996.

\bibitem{Kundur_D_1996_j-ieee-spm_blind_idr}
D.~Kundur and D.~Hatzinakos,
\newblock ``Blind image deconvolution revisited,''
\newblock {\em IEEE Signal Process. Mag.}, vol. 13, no. 6, pp. 61--63, Nov.
  1996.

\bibitem{Kato_M_1999_j-ieice-tfeccs_set-theoretic_bidbhsdm}
M.~Kato, I.~Yamada, and K.~Sakaniwa,
\newblock ``A set-theoretic blind image deconvolution based on hybrid steepest
  descent method,''
\newblock {\em IEICE Trans. Fund. Electron. Comm. Comput. Sci.}, vol. E82-A,
  no. 8, pp. 1443--1449, Aug. 1999.

\bibitem{Ahmed_A_2014_j-ieee-tit_blind_ducp}
A.~Ahmed, B.~Recht, and J.~Romberg,
\newblock ``Blind deconvolution using convex programming,''
\newblock {\em IEEE Trans. Inf. Theory}, vol. 60, no. 3, pp. 1711--1732, Mar.
  2014.

\bibitem{Campisi_P_2007_book_blind_idta}
P.~Campisi and K.~Egiazarian, Eds.,
\newblock {\em Blind Image Deconvolution: Theory and Applications},
\newblock CRC Press, 2007.

\bibitem{Comon_1996_j-ieee-spl_contrasts_mbd}
P.~Comon,
\newblock ``Contrasts for multichannel blind deconvolution,''
\newblock {\em Signal Process. Lett.}, vol. 3, no. 7, pp. 209--211, Jul. 1996.

\bibitem{Moreau_E_1997_j-ieee-spl_generlaized_cmbdls}
\'E. Moreau and J.-C. Pesquet,
\newblock ``Generalized contrasts for multichannel blind deconvolution of
  linear systems,''
\newblock {\em Signal Process. Lett.}, vol. 4, no. 6, pp. 182--183, Jun. 1997.

\bibitem{Zibulevsky_2001_j-neural-comput_blind_sssdsd}
M.~Zibulevsky and B.~A. Pearlmutter,
\newblock ``Blind source separation by sparse decomposition in a signal
  dictionary,''
\newblock {\em Neural Comput.}, vol. 13, no. 4, pp. 863--882, Apr. 2001.

\bibitem{Hoyer_P_2004_j-mach-learn-res_non-negative_mfsc}
P.~Hoyer,
\newblock ``Non-negative matrix factorization with sparseness constraints,''
\newblock {\em J. Mach. Learn. Res.}, vol. 5, pp. 1457--1469, 2004.

\bibitem{Hurley_N_2009_j-ieee-tit_comparing_ms}
N.~Hurley and S.~Rickard,
\newblock ``Comparing measures of sparsity,''
\newblock {\em IEEE Trans. Inf. Theory}, vol. 55, no. 10, pp. 4723--4741, Oct.
  2009.

\bibitem{Barak_B_2014_p-stoc_rounding_sofsr}
B.~Barak, J.~Kelner, and D.~Steurer,
\newblock ``Rounding sum-of-squares relaxations,''
\newblock in {\em Proc. ACM Symp. Theo. Comput. (STOC)}, New York, NY, USA, May
  31-Jun. 3, 2014.

\bibitem{Morup_M_2008_p-iscas_approximate_l0cnnmtf}
M.~M\o{}rup, K.~H. Madsen, and L.~K. Hansen,
\newblock ``Approximate {$L_0$} constrained non-negative matrix and tensor
  factorization,''
\newblock in {\em Proc. Int. Symp. Circuits Syst.}, May 2008, pp. 1328--1331.

\bibitem{Gray_W_1978_tr_variable_nd}
W.~C. Gray,
\newblock ``Variable norm deconvolution,''
\newblock Tech. {R}ep. SEP-14, Stanford Exploration Project, Apr. 1978,
\newblock \url{http://sepwww.stanford.edu/oldreports/sep14/14_19.pdf}.

\bibitem{Ji_H_2012_j-acha_image_ducsiwd}
H.~Ji, J.~Li, Z.~Shen, and K.~Wang,
\newblock ``Image deconvolution using a characterization of sharp images in
  wavelet domain,''
\newblock {\em Appl. Comp. Harm. Analysis}, vol. 32, no. 2, pp. 295--304, 2012.

\bibitem{Demanet_L_2014_j-information-inference_scaling_lrses}
L.~Demanet and P.~Hand,
\newblock ``Scaling law for recovering the sparsest element in a subspace,''
\newblock {\em Information and Inference}, 2014,
\newblock To appear.

\bibitem{Benichoux_A_2013_p-icassp_fundamental_pbdssip}
A.~Benichoux, E.~Vincent, and R.~Gribonval,
\newblock ``A fundamental pitfall in blind deconvolution with sparse and
  shift-invariant priors,''
\newblock in {\em Proc. Int. Conf. Acoust. Speech Signal Process.}, Vancouver,
  BC, Canada, May 26-31, 2013.

\bibitem{Krishnan_D_2011_p-cvpr_blind_dnsm}
D.~Krishnan, T.~Tay, and R.~Fergus,
\newblock ``Blind deconvolution using a normalized sparsity measure,''
\newblock in {\em Proc. IEEE Conf. Comput. Vis. Pattern Recogn.}, Colorado
  Springs, CO, USA, Jun. 21-25, 2011, pp. 233--240.

\bibitem{Esser_E_2013_j-siam-j-imaging-sci_method_fsssnnlspa}
E.~Esser, Y.~Lou, and J.~Xin,
\newblock ``A method for finding structured sparse solutions to non-negative
  least squares problems with applications,''
\newblock {\em SIAM J. Imaging Sci.}, vol. 6, no. 4, pp. 2010--2046, 2013.

\bibitem{Bolte_J_2013_j-math-prog-ser-a_proximal_almnnp}
J.~Bolte, S.~Sabach, and M.~Teboulle,
\newblock ``Proximal alternating linearized minimization fon nonconvex and
  nonsmooth problems,''
\newblock {\em Math. Progr. (Ser. A)}, Jul. 2013.

\bibitem{Chouzenoux_E_2013_tr_block_cvmfba}
E.~Chouzenoux, J.-C. Pesquet, and A.~Repetti,
\newblock ``A block coordinate variable metric forward-backward algorithm,''
\newblock Tech. {R}ep., 2013,
\newblock \url{http://www.optimization-online.org/DB_HTML/2013/12/4178.html}.

\bibitem{Sotthivirat_S_2002_j-ieee-tip_image_rupspscaa}
S.~Sotthivirat and J.~A. Fessler,
\newblock ``Image recovery using partitioned-separable paraboloidal surrogate
  coordinate ascent algorithms,''
\newblock {\em IEEE Trans. Image Process.}, vol. 11, no. 3, pp. 306--317, Mar.
  2002.

\bibitem{Chouzenoux_E_2014_j-optim-theory-appl_variable_mfbamsdfcf}
E.~Chouzenoux, J.-C. Pesquet, and A.~Repetti,
\newblock ``Variable metric forward-backward algorithm for minimizing the sum
  of a differentiable function and a convex function,''
\newblock {\em J. Optim. Theory Appl.}, vol. 162, no. 1, pp. 107--132, Jul.
  2014.

\bibitem{HiriartUrruty_J_1993_book_convex_ama}
J.-B. Hiriart-Urruty and C.~Lemar\'echal,
\newblock {\em Convex Analysis and Minimization Algorithms},
\newblock Springer-Verlag, 1993.

\bibitem{Combettes_P_2013_j-non-linear-anal_variable_mqfm}
P.~L. Combettes and B.~C. V\~u,
\newblock ``Variable metric quasi-{F}ej\'er monotonicity,''
\newblock {\em Nonlinear Anal.}, vol. 78, pp. 17--31, Feb. 2013.

\bibitem{Moreau_J_1965_j-bull-soc-math-fr_proximite_deh}
J.~J. Moreau,
\newblock ``Proximit\'e et dualit\'e dans un espace hilbertien,''
\newblock {\em Bull. Soc. Math. France}, vol. 93, pp. 273--299, 1965.

\bibitem{Chaux_C_2007_j-inv-prob_variational_ffbip}
C.~Chaux, P.~L. Combettes, J.-C. Pesquet, and V.~R. Wajs,
\newblock ``A variational formulation for frame based inverse problems,''
\newblock {\em Inverse Probl.}, vol. 23, no. 4, pp. 1495--1518, Aug. 2007.

\bibitem{Pustelnik_N_2011_j-ieee-tip_parallel_pairuhr}
N.~Pustelnik, C.~Chaux, and J.-C. Pesquet,
\newblock ``Parallel proximal algorithm for image restoration using hybrid
  regularization,''
\newblock {\em IEEE Trans. Image Process.}, vol. 20, no. 9, pp. 2450--2462,
  Sep. 2011.

\bibitem{Combettes_P_2011_incoll_proximal_smsp}
P.~L. Combettes and J.-C. Pesquet,
\newblock ``Proximal splitting methods in signal processing,''
\newblock in {\em Fixed-point algorithms for inverse problems in science and
  engineering}, H.~H. Bauschke, R.~Burachik, P.~L. Combettes, V.~Elser, D.~R.
  Luke, and H.~Wolkowicz, Eds., pp. 185--212. Springer Verlag, 2011.

\bibitem{Slavakis_K_2013_j-ieee-tsp_generalized_tosalus}
K.~Slavakis, Y.~Kopsinis, S.~Theodoridis, and S.~McLaughlin,
\newblock ``Generalized thresholding and online sparsity-aware learning in a
  union of subspaces,''
\newblock {\em IEEE Trans. Signal Process.}, vol. 61, no. 15, pp. 3760--3773,
  Aug. 2013.

\bibitem{Luo_Z_1992_j-optim-theory-appl_convergence_cdmcdm}
Z.~Q. Luo and P.~Tseng,
\newblock ``On the convergence of the coordinate descent method for convex
  differentiable minimization,''
\newblock {\em J. Optim. Theory Appl.}, vol. 72, no. 1, pp. 7--35, Jan. 1992.

\bibitem{Bolte_J_2010_p-icip_alternating_pabir}
J.~Bolte, P.~L. Combettes, and J.-C. Pesquet,
\newblock ``Alternating proximal algorithm for blind image recovery,''
\newblock in {\em Proc. Int. Conf. Image Process.}, Hong-Kong, China, Sep.
  26-29, 2010, pp. 1673--1676.

\bibitem{Xu_Y_2013_j-siam-j-imaging-sci_block_cdmrmoantfc}
Y.~Xu and W.~Yin,
\newblock ``A block coordinate descent method for regularized multiconvex
  optimization with applications to nonnegative tensor factorization and
  completion,''
\newblock {\em SIAM J. Imaging Sci.}, vol. 6, no. 3, pp. 1758--1789, 2013.

\bibitem{Allain_M_2006_j-ieee-tip_global_lchqa}
M.~Allain, J.~Idier, and Y.~Goussard,
\newblock ``On global and local convergence of half-quadratic algorithms,''
\newblock {\em IEEE Trans. Image Process.}, vol. 15, no. 5, pp. 1130--1142, May
  2006.

\bibitem{Osman_O_1996_book_seismic_ssem}
O.~S. Osman and E.~A. Robinson, Eds.,
\newblock {\em Seismic Source Signature Estimation and Measurement},
\newblock Number~18 in Geophysics Reprint Series. Soc. Expl. Geophysicists,
  Tulsa, OK, USA, 1996.

\bibitem{Ulrych_T_2005_book_information-based_ipa}
T.~Ulrych and M.~D. Sacchi,
\newblock {\em Information-based inversion and processing with applications},
\newblock Elsevier, 2005.

\bibitem{Walden_A_1986_j-geophys-prospect_nature_ngprcsd}
A.~T. Walden and J.~W.~J. Hosken,
\newblock ``The nature of the non-{G}aussianity of primary reflection
  coefficients and its significance for deconvolution,''
\newblock {\em Geophys. Prospect.}, vol. 34, no. 7, pp. 1038--1066, 1986.

\bibitem{Ricker_N_1940_j-geophysics_form_nswss}
N.~Ricker,
\newblock ``The form and nature of seismic waves and the structure of
  seismograms,''
\newblock {\em Geophysics}, vol. 5, no. 4, pp. 348--366, 1940.

\bibitem{Loganathan_P_2009_j-ieee-taslp_class_scaec}
P.~Loganathan, A.~W.~H. Khong, and P.~A. Naylor,
\newblock ``A class of sparseness-controlled algorithms for echo
  cancellation,''
\newblock {\em IEEE Trans. Audio Speech Lang. Process.}, vol. 17, no. 8, pp.
  1591--1601, Nov. 2009.

\bibitem{Drugman_T_2014_j-ieee-spl_maximum_pmslps}
T.~Drugman,
\newblock ``Maximum phase modeling for sparse linear prediction of speech,''
\newblock {\em Signal Process. Lett.}, vol. 21, no. 2, pp. 185--189, Feb. 2014.

\bibitem{Yukawa_M_2013_p-iswcs_sparsity-based_drpapfbsa}
M.~Yukawa, Y.~Tawara, S.~Sasaki, and I.~Yamada,
\newblock ``A sparsity-based design of regularization parameter for adaptive
  proximal forward-backward splitting algorithm,''
\newblock in {\em Proc. Int. Symp. Wireless Comm. Syst.}, Ilmenau, Germany,
  Aug. 27-30, 2013, pp. 1--4.

\bibitem{Chang_X_2014_PREPRINT_sparse_kmlilophddc}
X.~Chang, Y.~Wang, R.~Li, and Z.~Xu,
\newblock ``Sparse {K}-means with $\ell_{\infty}/\ell_0$ penalty for
  high-dimensional data clustering,''
\newblock {\em PREPRINT}, Mar. 2014.

\end{thebibliography}
\end{document}